\documentclass[12 pt]{article}

\usepackage{amsmath}
\usepackage{amsfonts}
\usepackage{amssymb}
\usepackage{amsthm}
\usepackage{graphicx}
\usepackage{xfrac}
\usepackage{tikz}

\begin{document}

\title{{\bf Genus theory and the factorization of class equations over $\mathbb{F}_p$}}
\author{Patrick Morton}
\date{January 30, 2019}
\maketitle

\begin{abstract} A new proof, depending only on genus theory, is given of a theorem of Stankewicz, which characterizes the primes $p$ for which the class equation $H_D(X)$ of the maximal order of the imaginary quadratic field $K=\mathbb{Q}(\sqrt{D})$ has a linear factor (mod $p$).  This yields a prime decomposition law for the primes $p$ with $p \nmid D$ in the real subfield of the Hilbert class field of $K$.
\end{abstract}

As is well-known, the Hilbert class equation (or class polynomial) is the polynomial $H_D(X)$ whose roots are the distinct $j$-invariants of elliptic curves with complex multiplication by the maximal order $R_K$ in the imaginary quadratic field $K=\mathbb{Q}(\sqrt{D})$ with discriminant $D$.  A root of $H_D(X)$ generates the Hilbert class field $\Sigma$ of $K$ over $K$.  The polynomial $H_D(X)$ always has a real root $\xi$, and over $\mathbb{Q}$ this root generates the real subfield $\Sigma_0=\mathbb{Q}(\xi)$ of $\Sigma$.  Recently, Stankewicz \cite{st} found a criterion for the polynomial $H_D(X)$ to have a root (mod $p$), for a given odd prime $p$ for which the Legendre symbol $(D/p)=-1$.  This criterion can be stated as follows.  \bigskip

\noindent {\bf Theorem} (Stankewicz).  {\it If $p$ is an odd prime for which $(D/p)=-1$, and $p$ does not divide the discriminant of $H_D(X)$, then $H_D(X)$ has a linear factor over $\mathbb{F}_p$ if and only if}
$$\left(\frac{-p}{q}\right) = 1, \ \ \forall \ q|D, \ q \ an \ odd \ prime.$$ \bigskip

Stankewicz derives this criterion from his analysis of rational $p$-adic points on twists of Shimura curves.  In this note I give a more direct proof of the criterion using genus theory and basic properties of the Hilbert class field.  The proof shows that the above theorem fits naturally into a discussion of genus theory.  (See \cite{co}, \cite{has2}, and \cite{ish}. Also see \cite{leo}, \cite{hak1}, \cite{hak2}.)  We will see that the same criterion holds for the prime $p=2$, if $2 \nmid D$.  The proof also yields a prime decomposition law for primes $p$, which do not divide $D$, in the real subfield of the Hilbert class field.

\section{Necessity.}

In this and the next section the integer $-N$ will denote the square-free part of the fundamental quadratic discriminant $D$, and $K$ is the imaginary quadratic field $K=\mathbb{Q}(\sqrt{-N})$ with discriminant $D$.  \bigskip

\noindent {\bf Theorem 1.} {\it Let $\Sigma_0$ denote the real subfield of the Hilbert class field $\Sigma$ of the quadratic field $K=\mathbb{Q}(\sqrt{-N})$.  Assume $p$ is an odd prime for which $\displaystyle \left(\frac{-N}{p} \right) =-1$.  If $p$ has a prime divisor $\mathfrak{p}$ of degree 1 in $\Sigma_0$, then}
$$\left(\frac{-p}{q}\right) = 1, \ \ \forall \ q|N, \ q \ an \ odd \ prime.$$

\noindent {\it This also holds for $p=2$ if $\displaystyle \left(\frac{-N}{p} \right) =-1$ is the Kronecker symbol, so that $D=-N \equiv 1$ (mod $4$)}. \bigskip

To prove this we use the decomposition
$$D = \prod_{q|D}{q^*}, \ q^* = (-1)^{(q-1)/2}q, \ q \ \textrm{odd}, \ 2^* = -4, 8, -8,$$
where the product is over all the prime divisors of $D$.  \medskip

The genus field of $K$ is the field $\Omega$, which is obtained by adjoining all the square-roots $\sqrt{q^*}$ to $K$, as $q$ varies over the prime divisors of $D$. It is the largest unramified extension of $K$ which is abelian over $\mathbb{Q}$, so that $\Omega \subseteq \Sigma$.  \medskip

Assume that the odd prime $p$ has a first degree prime divisor in $\Sigma_0$.  The conditions $\displaystyle \left(\frac{-N}{p} \right) =-1$ and $p$ odd imply that $p$ does not divide $D$, and $p$ has a first degree prime divisor in every subfield of $\Sigma_0$.  If $q$ is a prime $\equiv 1$ (mod $4$), then $q^*=q$, so $\mathbb{Q}(\sqrt{q}) \subseteq \Sigma_0$.  Hence, $p$ has a first degree prime divisor in $\mathbb{Q}(\sqrt{q})$, which implies that
$$\left(\frac{q}{p} \right) =1,  \ q \equiv 1 \ (\textrm{mod} \ 4), \ q|N.\eqno{(1.1)}$$
This implies then that
$$\left(\frac{-p}{q} \right) =1, \ q \equiv 1 \ (\textrm{mod} \ 4), \ q|N.\eqno{(1.2)}$$ \smallskip
If $2|D$ and $2^*=8$, the same argument also gives $\displaystyle \left(\frac{2}{p}\right)=1$, so $p \equiv \pm 1$ (mod $8$). \medskip

On the other hand, if there are several primes $q_i \equiv 3$ (mod $4$), $i = 1,2$, then $\sqrt{q_i^*} \in \Sigma$ implies that $\sqrt{q_1^* q_2^*} = \sqrt{q_1 q_2} \in \Sigma_0$.  Then $p$ has a first degree prime divisor in $\mathbb{Q}(\sqrt{q_1q_2})$, so we have
$$\left(\frac{q_1q_2}{p} \right) =1, \ q_i \equiv 3 \ (\textrm{mod} \ 4), \ q_i|N.$$ \smallskip
It follows that
$$\left(\frac{q_1^*}{p}\right) =\left(\frac{q_2^*}{p}\right), \ q_1 \equiv q_2 \equiv 3 \ (\textrm{mod} \ 4), \ q_i|N.\eqno{(1.3)}$$ \smallskip
We get a similar conclusion when $2|D$ and $2^* = -4, -8$, namely
$$\left(\frac{2^*}{p}\right) =\left(\frac{q^*}{p}\right), \ 2^*=-4,-8, \  q \equiv 3 \ (\textrm{mod} \ 4), \ q|N.\eqno{(1.4)}$$ \medskip

Now we use the fact that
$$\left(\frac{D}{p}\right)=\left(\frac{-N}{p}\right) = \prod_{q|D}{\left(\frac{q^*}{p}\right)}=-1.\eqno{(1.5)}$$ \smallskip
From (1.1), (1.3), (1.4) the terms with $q \equiv 1$ (mod $4$) or $q^*=8$ drop out, and we are left with
$$\left(\frac{q^*}{p}\right)^r = -1,$$\smallskip
where $r$ is the number of prime divisors of $D$ with $q\equiv 3$ (mod $4$) or $q^* = -4, -8$.  But this implies that $r$ is odd and $\displaystyle \left(\frac{q^*}{p}\right)=-1$ for all these prime divisors.  Hence,
$$\left(\frac{-p}{q}\right)=1, \ \textrm{if} \ q \equiv 3 \ (\textrm{mod} \ 4), \ q|N.\eqno{(1.6)}$$ \smallskip

If $p=2$, then $D$ is odd, and (1.2) and (1.3) still hold, where the symbols in (1.3) are now Kronecker symbols, while (1.4) falls away.  The same arguments show that (1.6) holds also in this case.  Together with (1.2), this proves Theorem 1. $\square$ \bigskip

\noindent {\bf Corollary 1.} {\it If the prime $p$ does not divide the discriminant of $H_D(X)$, $\displaystyle \left(\frac{-N}{p}\right)=-1$, and $H_{D}(X)$ (mod $p$) has a root in $\mathbb{F}_p$, then}
$$\left(\frac{-p}{q}\right) = 1, \ \ \forall \ q|N, \ q \ an \ odd \ prime.$$ \medskip

\noindent {\it Proof.} Since $p$ does not divide the discriminant of $H_D(X)$ and a real root of $H_D(X)$ generates $\Sigma_0$, it is clear that the factors of $H_D(X)$  mod $p$ correspond 1-1 to the prime divisors of $p$ in $\Sigma_0$.  The corollary is now immediate from Theorem 1. $\square$ \bigskip
 
\section{Sufficiency.}

Now we prove the converse of Theorem 1: \bigskip

\noindent {\bf Theorem 2.} {\it Let $\Sigma_0$ denote the real subfield of the Hilbert class field $\Sigma$ of the quadratic field $K=\mathbb{Q}(\sqrt{-N})$.  Assume $p$ is a prime for which $\displaystyle \left(\frac{-N}{p} \right) =-1$.  If $p$ satisfies the condition
$$\left(\frac{-p}{q}\right) = 1, \ \ \forall \ q|N, \ q \ \textrm{an odd prime},$$ \smallskip
then $p$ has a prime divisor $\mathfrak{p}$ of degree 1 in $\Sigma_0$.} $\square$ \bigskip

To prove this we consider the decomposition group of a prime divisor $\mathfrak{P}$ of $p$ in $\Sigma$.  First we note that if $\displaystyle \left(\frac{-p}{q}\right)=1$ for all odd prime divisors $q$ of $N$, then (1.1) holds, as does
$$\displaystyle \left(\frac{q^*}{p}\right)=-1, \ q\equiv 3 \ (\textrm{mod} \ 4), \ q|N, \ q \ \textrm{prime}.\eqno{(2.1)}$$ \smallskip
Now (1.5) implies that
$$\left(\frac{2^*}{p}\right)=(-1)^{r-1}\cdot(-1)=(-1)^r, \ \textrm{if} \ \ 2|D \ \textrm{and} \ p \neq 2,$$ \smallskip
where $r-1$ is the number of primes $q \equiv 3$ (mod $4$) dividing $N$.  But if $2|D$, then either: \medskip

$N \equiv 1$ (mod $4$), in which case $2^*=-4$ and $r-1$ is even, so that $r$ is odd, implying that $\displaystyle \left(\frac{2^*}{p}\right)=\left(\frac{-4}{p}\right)=-1;$ \medskip

$N \equiv 2$ (mod $8$), in which case $2^*=-8$ and $r$ is again odd, giving $\displaystyle \left(\frac{2^*}{p}\right)=\left(\frac{-8}{p}\right)=-1;$ \medskip

or $N \equiv 6$ (mod $8$), in which case $2^*=8$ and $r-1$ is odd, giving $\displaystyle \left(\frac{2^*}{p}\right)=\left(\frac{8}{p}\right)=1.$ \medskip

Thus, if $2|D$, we have (1.4) and the assertion in the sentence following (1.2).  This shows that $p$ splits completely in the real subfield $\Omega_0$ of the genus field $\Omega$.  (Note that $[\Omega: K]=2^{t-1}$, where $t$ is the number of distinct prime factors of $D$.  Thus $[\Omega_0:\mathbb{Q}]=2^{t-1}$, as well.)  Hence, the decomposition field of any prime divisor $\mathfrak{P}$ of $p$ in $\Sigma$ contains the field $\Omega_0$, and therefore the decomposition group $G_\mathfrak{P}$ is contained in $H=\textrm{Gal}(\Sigma/\Omega_0)$.  \medskip

It suffices to show $G_\mathfrak{P}=\{1, \tau \}$ for some $\mathfrak{P}|p$, where $\tau$ is complex conjugation.  If this holds, then $\Sigma_0$, which is the fixed field of $\tau$, is the largest field in which the prime below $\mathfrak{P}$ has degree 1, i.e. $\mathfrak{p}=\mathfrak{P} \mathfrak{P}^\tau$ is a first degree prime divisor of $p$ in $\Sigma_0$.  \medskip

Let $J = H \ \cap \ \textrm{Gal}(\Sigma/K)$ be the subgroup of $\textrm{Gal}(\Sigma/K)$ corresponding to $\Omega$ in the Galois correspondence, so that $[H:J]=2$ and $H = J \cup J\tau$.  By the genus theory \cite{has2}, $J$ corresponds to the subgroup of squares in $Pic(R_K)$, in the Artin correpondence between ideal classes in $R_K$ and elements of the Galois group $\textrm{Gal}(\Sigma/K)$.  \medskip

We now have what we need to complete the proof.  The prime $p$ is inert in $K$, so that $(p)=pR_K$ is a principal ideal and therefore $p$ splits completely in the extension $\Sigma/K$.  In particular, the prime divisors $\mathfrak{P}$ of $p$ in $\Sigma$ have degree $2$ over $\mathbb{Q}$, so that $|G_\mathfrak{P}|=2$.  Furthermore, we know that $K$ is not contained in the decomposition field of any $\mathfrak{P}$, and therefore $G_\mathfrak{P} \not \subset J \subset \textrm{Gal}(\Sigma/K)$.  Hence, $G_\mathfrak{P} \subset H$ is generated by some $\sigma \tau$, with $\sigma \in J$.  But $\sigma = \psi^{2}$ for some $\psi \in \textrm{Gal}(\Sigma/K)$, by the characterization of the group $J$, and $\psi^{-1} G_\mathfrak{P} \psi = G_{\mathfrak{P}^{\psi}}=\{1, \psi^{-1} \sigma \tau \psi \}$, with $\psi^{-1} \sigma \tau \psi=\psi^{-2} \sigma \tau = \tau$.  This shows that $G_{\mathfrak{P}^{\psi}} = \{1, \tau \}$ and completes the proof. $\square$ \bigskip

\noindent {\bf Corollary 2.} {\it If $p$ satisfies $\displaystyle \left(\frac{-N}{p}\right)=-1$ and $\displaystyle \left(\frac{-p}{q}\right) = 1$, for all odd primes $q$ such that $q|N$, then $H_{D}(X)$ (mod $p$) has a root in $\mathbb{F}_p$.} \medskip

\noindent {\it Proof.} Theorem 2 implies that $H_D(X)$ has a linear factor over $\mathbb{Q}_p$ and therefore $H_D(X)$ (mod $p$) has a root in $\mathbb{F}_p$. $\square$ \bigskip

Note that the hypothesis $\displaystyle  \left(\frac{-N}{p}\right)=-1$ of Corollary 2 holds for all the prime divisors of the discriminant of $H_D(X)$ which do not divide $D$, by a result of Deuring \cite{d}.

\section{A prime decomposition law for $\Sigma_0$.}

The proof of Theorem 2 shows that when $p$ has a first degree prime divisor $\mathfrak{p}$ in $\Sigma_0$, then it has as many prime divisors of degree $1$ as there are distinct elements $\psi$ in $\textrm{Gal}(\Sigma/K)$ for which $ G_{\mathfrak{P}^{\psi}} = \psi^{-1} G_\mathfrak{P} \psi = \{1, \tau\}$, where $\mathfrak{P}|\mathfrak{p}$.  (Note that the $\psi$ are in distinct cosets of $G_{\mathfrak{P}}= \{1, \tau\}$, so the prime divisors $\mathfrak{P}^{\psi}$ are distinct.)  This holds if and only if $\psi^{-1} \tau \psi = \tau$, i.e., if and only if $\psi^2=1$.  The number of such elements $\psi$ is exactly $2^{t-1}$, since this is the order of the $2$-Sylow subgroup of the class group $Pic(R_K)$.  Thus we have: \bigskip

\noindent {\bf Theorem 3.} {\it If $p$ is a prime for which $\displaystyle \left(\frac{-N}{p}\right) = -1$, and $p$ has a prime divisor of degree $1$ in $\Sigma_0$, then it has exactly $2^{t-1}$ such prime divisors, where $t$ is the number of distinct prime factors of $D$.} $\square$ \bigskip

Taken together, Theorems 1-3 yield the following decomposition law for the real subfield $\Sigma_0$ of $\Sigma$.  \bigskip

\noindent {\bf Prime Decomposition Law in $\Sigma_0$.} {\it Let $p$ be a prime that does not divide the discriminant $D$ of the imaginary quadratic field $K=\mathbb{Q}(\sqrt{D})$. \medskip

(a) If $\displaystyle \left(\frac{D}{p}\right)=1$, then in $\Sigma_0$, $p$ splits into $h(D)/f$ primes of degree $f$ over $\mathbb{Q}$, where $f$ is the order of a prime ideal divisor $\wp$ of $p$ in $Pic(R_K)$. \medskip

(b) If $\displaystyle \left(\frac{D}{p}\right)=-1$ and $\displaystyle \left(\frac{-p}{q}\right)=1$ for all odd prime divisors $q$ of $D$, then in $\Sigma_0$, $p$ splits into $r_1=2^{t-1}$ primes of degree 1 and $r_2=(h(D)-2^{t-1})/2$ primes of degree 2 over $\mathbb{Q}$. \medskip

(c) If $\displaystyle \left(\frac{D}{p}\right)=-1$ and $\displaystyle \left(\frac{-p}{q}\right)=-1$ for some odd prime divisor $q$ of $D$, then in $\Sigma_0$, $p$ splits into $h(D)/2$ primes of degree 2 over $\mathbb{Q}$. } $\square$ \bigskip

This law immediately implies the following density result.  \bigskip

\noindent {\bf Theorem 4.} {\it The density of primes $p \in \mathbb{Z}^+$ for which $H_D(X)$ has a linear factor (mod $p$) is $\displaystyle d(P(\Sigma_0)) = \frac{1}{2h(D)}+\frac{1}{2^t}$, where $t$ is the number of distinct prime factors of $D$.}  \bigskip

\section{Discriminant divisors.}

The Prime Decomposition Law proved in the last section raises the question: how do the prime divisors $p$ of the discriminant $D$ split in the real subfield $\Sigma_0$? \medskip

While the Prime Decomposition Law in Section 2 refers to the quadratic character of $p$ with respect to the prime factors of the discriminant, a similar law for the prime factors $p$ of $D$ cannot simply involve these same prime factors.  For example, we have the following congruence for the class equation of discriminant $D=-8p$, for $p>13$, from \cite{mor2}, Theorem 1.1:
\begin{align*}
H_{-8p}(t) & \equiv (t-1728)^{2\epsilon_1}(t-8000)^{2\epsilon_2}(t+3375)^{4\epsilon_3}H_{-15}(t)^{4\epsilon_4}\\
& \times \prod_{i}{(t^2+a_i t+b_i)^2} \ (\textrm{mod} \ p).
\end{align*}
The product is over certain distinct, irreducible quadratic factors, different from $H_{-15}(t)=t^2+191025t-121287375$ (mod $p$).  The exponents in the congruence for $H_{-8p}(t)$ are given by
\begin{align*}
\epsilon_1&=\frac{1}{2} \left(1-\left(\frac{-4}{p}\right)\right), \ \ \epsilon_2=\frac{1}{2} \left(1-\left(\frac{-8}{p}\right)\right),\\
\epsilon_3&=\frac{1}{2}\left(1-\left(\frac{-7}{p}\right)\right), \ \ \epsilon_4=\frac{1}{4}\left(1-\left(\frac{-15}{p}\right)\right) \left(1-\left(\frac{5}{p}\right)\right).
\end{align*}
Thus, $H_{-8p}(t)$ has a linear factor modulo $p$ if and only if one of the exponents $\epsilon_1, \epsilon_2, \epsilon_3$ is $1$, which is the case (for $p>13$) exactly when
$$p \equiv 3, 5, 7 \ (\textrm{mod} \ 8) \ \ \textrm{or} \ \ \left(\frac{-p}{7}\right)=+1, \ \textrm{i.e.}, \ p \equiv 3, 5, 6 \ (\textrm{mod} \ 7).$$
Thus, the quadratic character of $p$ with respect to the prime $q=7$ is also relevant in this case! \medskip

Similar congruences are proven in \cite{mor2} for the discriminants $D = -3p$ (when $p \equiv 1$ (mod $4$)) and $D=-12p$ (when $p \equiv 3$ (mod $4$)).  If $p>53$ is a prime with $p \equiv 1$ (mod 4), then we have the following congruences (mod $p$):
$$H_{-3p}(t)\equiv (t-54000)^{2\delta_1}(t-8000)^{4\delta_2}H_{-20}(t)^{2\delta_4}H_{-32}(t)^{2\delta_5}\prod_{i}(t^2+c_i t+d_i)^2,$$
$$H_{-12p}(t)\equiv t^{2\delta_1}(t+32768)^{4\delta_3}H_{-20}(t)^{2\delta_4}H_{-32}(t)^{2\delta_5}H_{-35}(t)^{4\delta_6}\prod_{j}{(t^2+c_j t+d_j)^2}.$$
Here the quadratic polynomials $t^2+c_i t+d_i$ and $t^2+c_j t+d_j$ and class equations $H_{-20}(t)$, $H_{-32}(t)$, and $H_{-35}(t)$ (also quadratic) are distinct and irreducible (mod $p$) when they occur, and the exponents on the linear factors are given by
$$\delta_1=\frac{1}{2}\left(1-\left(\frac{-3}{p}\right)\right), \ \ \delta_2=\frac{1}{2}\left(1-\left(\frac{-8}{p}\right)\right), \ \ \delta_3=\frac{1}{2}\left(1-\left(\frac{-11}{p}\right)\right).$$
Thus, whether $H_{-3p}(t)$ has a linear factor (mod $p$) in this case depends on $\left(\frac{-8}{p}\right)=\left(\frac{p}{2}\right)$, i.e., the residue class of $p$ (mod $8$), even though the discriminant $D=-3p$ is odd.  Further, whether $H_{-12p}(t)$ (the class equation for a non-maximal order) has a linear factor (mod $p$) depends on $\left(\frac{-11}{p}\right)=\left(\frac{p}{11}\right)$.  When $p \equiv 3$ (mod $4$) and $p >53$, then in place of the two congruences above, we have the single congruence
\begin{align*}
H_{-12p}(t) & \equiv t^{2\delta_1}(t-54000)^{2\delta_1}(t-8000)^{4\delta_2}(t+32768)^{4\delta_3} H_{-20}(t)^{4\delta_4}H_{-32}(t)^{4\delta_5}\\
& \times H_{-35}(t)^{4\delta_6} \prod_{i}{(t^2+c_i t+d_i)^2} \ (\textrm{mod} \ p),
\end{align*}
and so whether $H_{-12p}(t)$ has a linear factor (mod $p$) depends again on the prime $q=11$. \medskip

Also note that the class equations $H_{-p}(X)$ and $H_{-4p}(X)$ always have linear factors modulo $p$, by results of Elkies and Kaneko.  See \cite{brm}, p. 95. \medskip

In general, if $D=-pQ$, is there a criterion involving the primes dividing $Q(Q-1)$ which determines when $H_D(X)$ has a linear factor (mod $p$)?

\noindent Dept. of Mathematical Sciences \smallskip

\noindent Indiana University - Purdue University at Indianapolis (IUPUI) \smallskip

\noindent 402 N. Blackford St., Indianapolis, IN, USA 46202 \smallskip

\noindent \em{e-mail: pmorton@math.iupui.edu}


\bigskip

\begin {thebibliography}{WWW}

\bibitem[1]{brm} J. Brillhart and P. Morton, Class numbers of quadratic fields, Hasse invariants of elliptic curves, and the supersingular polynomial, J. Number Theory 106 (2004), 79-111.

\bibitem[2]{co} D. A. Cox, {\it Primes of the form $x^2+ny^2$}, John Wiley and Sons, New York, 1989.

\bibitem[3]{d} M. Deuring, Teilbarkeitseigenschaften der singul\"aren Moduln der elliptischen Funktionen und die Diskriminante der Klassengleichung, Commentarii Mathematici Helvetici 19 (1946), 74-82.

\bibitem[4]{hak1} F. Halter-Koch, Arithmetische Theorie der Normalk\"orper von 2-Potenzgrad mit Diedergruppe, J. Number Theory 3 (1971), 412-443.

\bibitem[5]{hak2} F. Halter-Koch, Geschlechtertheorie der Ringklassenk\"orper, J. reine angew. Math. 250 (1971), 107-108.

\bibitem[6]{has2} H. Hasse, Zur Geschlechtertheorie in quadratischen Zahlk\"orpern, J. Math. Soc. Japan 3 (1951), 45-51.

\bibitem[7]{ish} M. Ishida, The Genus Fields of Algebraic Number Fields, Lecture Notes in Mathematics 555, Springer, Berlin, 1976.

\bibitem[8]{leo} Leopoldt, H., Zur Geschlechtertheorie in abelschen Zahlk\"orpern, {\it Math. Nachr.} {\bf 9} (1953), 351-362.

\bibitem[9]{mor2} P. Morton, Explicit congruences for class equations, Functiones et Approximatio Commentarii Math. 51 (2014), 77-110.

\bibitem[10]{st} J. Stankewicz, Twists of Shimura curves, Canad. J. Math. 66 (2014), 924-960; and at http://arxiv.org/abs/1208.3594.

\end{thebibliography}
\end{document}